\documentclass[11pt]{article}
\usepackage{graphicx}
\usepackage{amsmath}
\usepackage{amsfonts}
\usepackage{symbols}
\usepackage{booktabs}
\usepackage[margin=20truemm]{geometry}
\graphicspath{{./figures/}}
\title{The Galerkin method for a regularised combined field integral equation without a dual basis function}
\author{Kazuki~Niino, Shunpei~Yamamoto}
\begin{document}
\maketitle
\begin{abstract}
We propose discretisation of a regularised combined field integral equation~(regularised CFIE) only with the Rao-Wilton-Glisson~(RWG) basis function.
The CFIE is a formulation of integral equations, which avoids the so-called ficticious frequencies of integral equations.
The most typical CFIE, which is a linear combination of the electric field integral equation~(EFIE) and magnetic field integral equation~(MFIE),
is known to be ill-conditioned and requires many iterations when solved with iteration methods such as the generalised minimum residual~(GMRES) method.
The regularised CFIE is another formulation of the CFIE to solve this problem
by applying a regularising operator to the part of the EFIE.
In several previous studies the regularising operator is determined based on the Calder\'{o}n preconditioning.
This regularising operator however takes much more computatonal time than the standard CFIE
since discretising the EFIE with the Calder\'{o}n preconditioner requires the dual basis function.
In this article we propose a formulation of the regularised CFIE, which can be discretised with the Galerkin method without the dual basis function.
\end{abstract}

\section{Introduction}
Numerical methods based on boundary integral equations are effective for solving electromagnetic wave scattering problems.
The electric field integral equation~(EFIE) and magnetic field integral equation~(MFIE), which are ones of simplest formulations of boundary integral equations for Maxwell's equations,
are konwn to suffer from the so-called ficticious frequencies, which correspond to resonance frequencies inside a domain.
The combined field integral equation~(CFIE) is one of formulations of integral equations free from the ficticious frequencies.
The CFIE is typically obtained as a linear combination of the EFIE and MFIE, and indeed has no real ficticious frequencies.
However the CFIE is usually ill-conditioned due to the hyper-singular operator in the EFIE,
and iterative linear solvers require many iteration numbers to solve the CFIE.
In order to solve this problem the regularised CFIE, which regularises the hyper-singular operator by applying an appropriate operator (called regularising operator) to the counterpart of EFIE,
has been studied \cite{adams2004,bagci2009,bruno2009,boubendir2014}.
Among several options of the regularising operators, one based on the Calder\'{o}n formula \cite{nedelec2001} is a promising regularising operator.
This method utilises the electric field integral operator~(EFIO) as the regularising operator and multiply it to the part of the EFIE
since the Calder\'{o}n formula implies that the square of EFIO is well-conditioned.
The regularisation based on the Calder\'{o}n formula has been widely studied for preconditioning other integral equations such as the EFIE \cite{christiansen2003,andriulli2008} and PMCHWT \cite{yan2010,cools2011}.
This method is also effective for regularising the CFIE and successfully reduces the number of iterations.
However the computational cost of a matrix corresponding to discretisation of the regularising operator is much more than that of the standard EFIE. 
This fact is true for both of the two typical discretisation method of integral equations, namely the Galerkin method and Nystr\"{o}m method.
In the Galerkin method one has to use two mutually-dual basis functions (e.g.\ the Rao-Wilton-Glisson~(RWG) basis function \cite{rao1982} and Buffa-Christiansen~(BC) \cite{buffa2007} basis function) to discretise the square of the EFIO.
The BC basis function is defined as a linear combination of the RWG on the barycentric reinfement of the orignal mesh \cite{andriulli2008}.
Hence the operator discretised with the BC basis function requires much more computational time since the barycentric reinement in general has 6 times more number of elements than the original mesh.
In the Nystr\"{o}m method, the computaion of the squared EFIO is also difficult due to its strong singularity although several techniques to compute it has been proposed \cite{contopanagos2002}.

Another regularising operator for CFIE is the scalar potential, the weak singular part of the EFIO.
In \cite{bruno2009} this method was proposed in order to avoid the difficulty of the Nystr\"{o}m method stated in the previous paragraph.
Discretising this regularising operator with the Nystrom method is easier than that based on the Calder\'{o}n formula since this regularising operator no longer contains the hyper-singular operator.
In this article we propose that this regularising operator is also effective with the Galerkin method.
We will show that the regularised CFIE with this operator can be discretised with only the RWG function, which leads to reduction of the computational time and memory requirement.

\setlength{\tabcolsep}{4pt} 
\section{formulations of CFIE}
We start with making formulation of CFIEs.
A typical formulation of CFIE is given as follows:
\begin{align}\label{eq:old_cfie}
\left\{\left(\frac{\mathcal{I}}{2} - \mathcal{K}_k\right) + \zi\omega\mu\alpha \bn\times \mathcal{T}_k\right\}\bj = \bn\times\bH^{\mathrm{inc}} + \alpha \bn\times(\bE^{\mathrm{inc}}\times\bn),
\end{align}
which is the sum of EFIE
\begin{align*}
  \zi\omega\mu\mathcal{T}_k\bj = \bE^{\mathrm{inc}}\times\bn
\end{align*}
and MFIE
\begin{align*}
  \left(\frac{\mathcal{I}}{2} - \mathcal{K}_k\right)\bj = \bn\times\bH^{\mathrm{inc}}
\end{align*}
where
\begin{align*}
\mathcal{T}_k\bj &=  \left(\mathcal{T}^S_k + \frac{1}{k^2}\mathcal{T}^N_k\right)\bj\\
\mathcal{T}^S_k\bj &= \bn\times\int_\Gamma G^k(x,y)\bj(y)\md S_y,\\
\mathcal{T}^N_k\bj &= \bn\times\int_\Gamma {\nabla\nabla}G^k(x,y)\bj(y)\md S_y,\\
\mathcal{K}_k\bm &= \bn\times\int_\Gamma \nabla^x G^k(x,y)\times\bm(y)\md S_y,
\end{align*}
$\bj$ is the unknown surface current, $\Gamma$ is the surface of a scatterer, $\bn$ is the outward unit normal vector on $\Gamma$, $k$ is the wave number, 
$\bE^{\mathrm{inc}}$ and $\bH^{\mathrm{inc}}$ are respectively the electric and magnetic field of the incident wave,
and $G^k$ is the fundamental solution of the Helmholtz equation in 3D:
\begin{align*}
  G^k(x,y) = \frac{\ex^{ik|x-y|}}{4\pi|x-y|}.
\end{align*}
It is known that the CFIE avoids the ficticious frequencies with $\operatorname{Re}\alpha\neq 0$.
The CFIE in \eqref{eq:old_cfie} is however ill-conditioned, and the simple discretisation of \eqref{eq:old_cfie} results in many iteration numbers when solved with a Krylov subspace method.
As a solution to this problem, the regularised CFIE, in which the term $\alpha\bn\times\cdot$ is replaced with an appropriate regularising operator $\mathcal{R}$, namely,
\begin{align}\label{eq:reg_cfie}
\left\{\left(\frac{\mathcal{I}}{2} - \mathcal{K}_k\right) + \zi\omega\mu\mathcal{R} \mathcal{T}_k\right\}\bj = \bn\times\bH^{\mathrm{inc}} + \mathcal{R}(\bE^{\mathrm{inc}}\times\bn)
\end{align}
has been studied \cite{bagci2009,bruno2009,boubendir2014}.
In what follows we introduce some choises of the regularising operator $\mathcal{R}$ and discretisation of the regularised CFIE based on the Galerkin method.

\section{Discretisation of the regularised CFIE based on Galerkin's method}
\subsection{Conventional regularised CFIE}\label{sec:conv_reg_cfie}
\begin{figure}[htbp]
  \begin{center}
    \includegraphics[width=50mm]{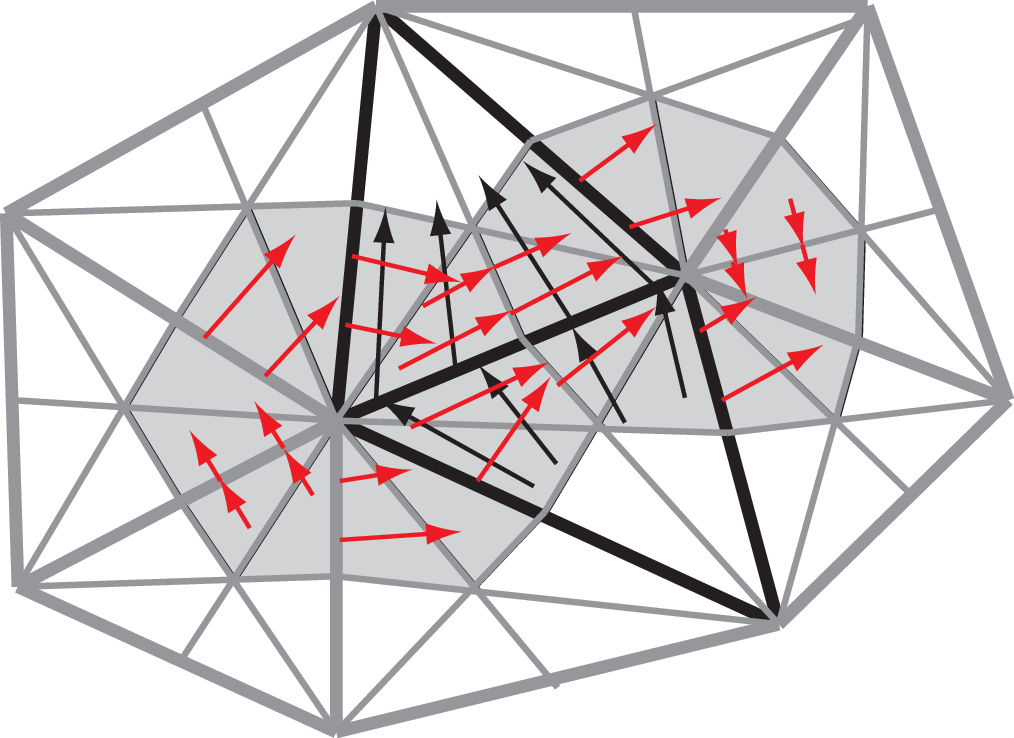}
    \caption{The RWG and BC basis functions.
      The thick lines forms a triangular mesh and thin lines are its barycentric refinement. The BC function (red arrows) is defined as a linear combination of the RWG on the refined mesh.
      The black arrow is the RWG function defined on the original mesh.}
    \label{figure:bc_basis}
  \end{center}
\end{figure}

In many previous studies the regularising operator $\mathcal{R}$ is decided based on the Calder\'{o}n formula \cite{nedelec2001}:
\begin{align}\label{eq:calderon}
  k^2\mathcal{T}_k^2 = \frac{\mathcal{I}}{4} - \mathcal{K}_k^2.
\end{align}
This equation implies that the square of the operator $\mathcal{T}_k$ is proportional to the identity operator up to compact purterbation since the operator $\mathcal{K}_k$ is compact if the surface $\Gamma$ is smooth.
Hence choosing the regularising operator $\mathcal{R}$ as a compact purterbation of $\mathcal{T}_k$ is a promising option.
However the simplest regularising operator in this sence $\mathcal{R}=\mathcal{T}_k$ makes the regularised CFIE in \eqref{eq:reg_cfie} suffer from the fictitious frequency again.
Indeed, by substituting $\mathcal{R}_k=\mathcal{T}_k$ into \eqref{eq:reg_cfie} and applying the Calder\'{o}n formula in \eqref{eq:calderon},
the LHS of the regularised CFIE in \eqref{eq:reg_cfie} leads to
\begin{align*}
  \left(\frac{\mathcal{I}}{2} + \mathcal{K}_k\right)^2\bj,
\end{align*}
which is sqaured magnetic field integral operator~(MFIO).
Consequently previous studies utilises compact purterbations of $\mathcal{T}_k$ as the regularising operator 
such as $\mathcal{R}=\omega\varepsilon\mathcal{T}_{ik}$:
\begin{align}\label{eq:reg_cfie_conv}
\left\{\left(\frac{\mathcal{I}}{2} - \mathcal{K}_k\right) + \zi k^2\mathcal{T}_{ik} \mathcal{T}_k\right\}\bj = \bn\times\bH^{\mathrm{inc}} + \omega\varepsilon\mathcal{T}_{ik}(\bE^{\mathrm{inc}}\times\bn).
\end{align}
This integral equation can be discretised as follows:
\begin{align}\label{eq:disc_reg_cfie_conv}
  (\widetilde{K} + \zi k^2\widetilde{R} \widetilde{G}^{-1} T  ) \bx = \widetilde{\bb},
\end{align}
where
\begin{align}
  \nonumber (\widetilde{K})_{ij} &= \left\langle \bn\times\bss_i , \left(\frac{\mathcal{I}}{2}-\mathcal{K}_k\right)\bt_j \right\rangle,\\
  \nonumber (\widetilde{R})_{ij} &= \langle \bn\times\bss_i, \mathcal{T}_{ik}\bss_j\rangle,\\
  \nonumber (\widetilde{G})_{ij} &= \langle \bn\times\bt_i, \bss_j\rangle,\\
  \label{eq:operator_t} (T)_{ij} &= \langle \bn\times\bt_i, \mathcal{T}_{k}\bt_j\rangle,\\
  \nonumber (\widetilde{\bb})_i &= \langle \bn\times\bss_i, \bn\times\bH^{\mathrm{inc}} + \omega\varepsilon\mathcal{T}_{ik}(\bE^{\mathrm{inc}}\times\bn) \rangle
\end{align}
and $\langle \cdot , \cdot \rangle$ is the $L^2(\Gamma)$ inner product defined as
\begin{align*}
  \langle \bu, \bv \rangle &= \int_\Gamma \overline{\bu}\cdot \bv \md S.
\end{align*}
The functions $\bt_i$ and $\bss_i$ are basis functions used as testing and trial functions in each matrix.
The most naive choise of the basis functions, namely a single basis funcions such as the RWG function for both $\bt_i$ and $\bss_i$, is known to be breakdown
since the Gram matrix $\widetilde{G}$ is singular and cannot be inverted in this case.
Hence we usually need two sets of mutually dual basis functions such as the RWG function for $\bt_i$ and the BC basis functions for $\bss_i$ defined as in Fig.~\ref{figure:bc_basis}.
Since the BC basis function is defined as a linear combination of the RWG function on the barycentric refinement,
a matrx including the BC function is usually computed with a matrix which is discretisation of an operator on the barycentric refinement,
and transformation between the orignal mesh and its barycentric refinement.
For example the matrix $\widetilde{R}$ is calculated with
\begin{align*}
  \widetilde{R} = C^T R_b C
\end{align*}
where
\begin{align*}
  (R_b)_{ij} &= \langle \bn\times\bt^b_i, \mathcal{T}_{ik} \bt^b_j\rangle,\\
\end{align*}
$\bt^b_i$ is the RWG function defined on the barycentric refinement and $C$ is the transformation matrix between $\bss_i$ and $\bt^b_i$ \cite{andriulli2008}.
The size of the matrix $R_b$ is approximately 6 times more than the sizes of the matrices $\widetilde{R}, T$ and $\widetilde{K}$ since the barycentric refinement has 6 times more triangular elements of the original mesh.
Consequently, in order to compute the matrix $\widetilde{R}$,  one has to deal with a matrix with larger size.
This leads to much more computainal cost for calculating the coefficient matrix in \eqref{eq:disc_reg_cfie_conv} than that of the standard CFIE in \eqref{eq:old_cfie}
although the coefficient matrix is well-conditioned.

\subsection{Proposed method}
In this section we propose a formulation of the regularised CFIE, which can be discretised only with the RWG basis function.
The operator $\mathcal{R}=\omega\varepsilon\mathcal{T}^S_0$ is utilised as the regularising operator in this formulation:
\begin{align}\label{eq:reg_cfie_new}
\left\{\left(\frac{\mathcal{I}}{2} - \mathcal{K}_k\right) + \zi k^2\mathcal{T}^S_{0} \mathcal{T}_k\right\}\bj = \bn\times\bH^{\mathrm{inc}} + \omega\varepsilon\mathcal{T}^S_{0}(\bE^{\mathrm{inc}}\times\bn)
\end{align}
The regularising operator $\mathcal{R}=\omega\varepsilon\mathcal{T}^S_0$ was proposed in \cite{bruno2009} in order to discretise the regularised CFIE with the Nystr\"{o}m method,
which is one of the collocation methods to make use of integral points of a quadrature rule as the collocation points.
Because the Nystr\"{o}m method calculates integrals with the quadrature rule, dealing with hyper-singular integrals, in particular the square of the operator $\mathcal{T}^N_k$, which appears in the conventional regularised CFIE in \eqref{eq:reg_cfie_conv}, is difficult.
Although some remedies to compute those hyper-singular operators have been proposed (such as in \cite{contopanagos2002}), they require more computational cost.
Hence the regularised CFIE in \eqref{eq:reg_cfie_new} with the regularising operator $\mathcal{R}=\mathcal{T}^S_0$ was proposed in \cite{bruno2009} as a better regularising operator for the Nystr\"{o}m method.
In the Galerkin method this regularising operator can be exploited in order to discretise the integral equation with only the RWG basis function.
In fact the regularised CFIE in \eqref{eq:reg_cfie_new} is discretised as follows:
\begin{align}\label{eq:disc_reg_cfie_new}
(K + \zi k^2RG^{-1}T)\bx = \bb
\end{align}
where
\begin{align*}
(K)_{ij} &= \left\langle \bt_i, \left(\frac{\mathcal{I}}{2} - \mathcal{K}_k\right)\bt_j\right\rangle,\\
(R)_{ij} &= \langle \bt_i, \mathcal{T}^S_0(\bn\times\bt_j)\rangle,\\
(G)_{ij} &= \langle \bn\times\bt_i, \bn\times\bt_j\rangle = \langle \bt_i, \bt_j\rangle,\\
(\bb)_i &= \langle \bt_i, \bn\times\bH^{\mathrm{inc}} + \omega\varepsilon\mathcal{T}^S_{0}(\bE^{\mathrm{inc}}\times\bn)\rangle,
\end{align*}
and $T$ is defined in \eqref{eq:operator_t}.
In this formulation the descritisation of the regularising operator $\mathcal{T}^S_0$ in matrix $R$ is peculiar
while the matrices $T$ and $K$ are stardard discretisation of the EFIO and MFIO.
We adopt the discretisation of the regularising operator in $R$ so that the Gram matrix $G$ is regular without the use of the dual basis function.
This discretisation can be achieved thanks to the replacement of the regularising operator from $\mathcal{R}=\mathcal{T}_{ik}$ to $\mathcal{R}=\mathcal{T}^S_0$.
In fact the operator $\mathcal{T}_{ik}$ discretised in this way
\begin{align*}
  \langle \bt_i, \mathcal{T}_{ik}(\bn\times\bt_j)\rangle,
\end{align*}
diverges due to the hyper singular operator $\mathcal{T}_{ik}^N$ contained in $\mathcal{T}_{ik}$.
The matrix $R$ is interpreted as the discretisation of $\mathcal{T}^S_0$ regarding it as an operator from $H_{\mathrm{curl}}^{-\frac12}$ to $H_{\mathrm{curl}}^{-\frac12}$,
and $\mathcal{T}^S_0$ is indeed such an operator while $\mathcal{T}_{ik}$ is not due to the hyper-singular operator $\mathcal{T}^N_{ik}$.
In this way the matrix $R$ can be calculated without divergence.
In this discretisation all the matrices contain only the single basis function $\bt_i$ and the dual basis function is not necessarry.
Thus computation of larger matrices such the matrix $R_b$ in the conventinal method introduced in section \ref{sec:conv_reg_cfie} is not required.
Hence the computaional cost for solving \eqref{eq:disc_reg_cfie_new} is expected to be less than that for \eqref{eq:disc_reg_cfie_conv}.
Also it is prooved that the operator
\begin{align*}
 \left(\frac{\mathcal{I}}{2} + \mathcal{K}_k\right) + \zi k^2\mathcal{T}^S_0 \mathcal{T}_k,
\end{align*}
which is the operator in the LHS of \eqref{eq:reg_cfie_new}, is bounded with compact perturbation\cite{bruno2009}.

\begin{table*}[!htb]
  \caption{Relative errors of numerical methods.}
  \label{table:rel_err}
  \centering
  \begin{tabular}{ccllllllll}
    \toprule
    && \multicolumn{2}{l}{eq.\ \eqref{eq:disc_reg_cfie_new}, $N=5780$} && \multicolumn{2}{l}{eq.\ \eqref{eq:disc_reg_cfie_new}, $N=30420$} && \multicolumn{2}{l}{eq.\ \eqref{eq:disc_reg_cfie_conv}, $N=5780$} \\
    \cmidrule(lr){3-4}    \cmidrule(lr){6-7}    \cmidrule(lr){9-10}
    wave number-$k$ && It. & Error && It. & Error && It. & Error \\
    \midrule
    1.0 && 11 & 1.503\% &  & 11 & 0.6546\% & & 13 & 1.493\%\\
    2.75&& 13 & 3.982\% &  & 13 & 1.736\% & & 15 & 3.982\%\\
    6.1 && 16 & 8.092\% &  & 16 & 3.530\% & & 18 & 8.119\%\\
    \bottomrule
  \end{tabular}
\end{table*}
\section{Numerical example}\label{sec:numerical_example}
\begin{figure}[htbp]
  \begin{center}
    \includegraphics[width=90mm]{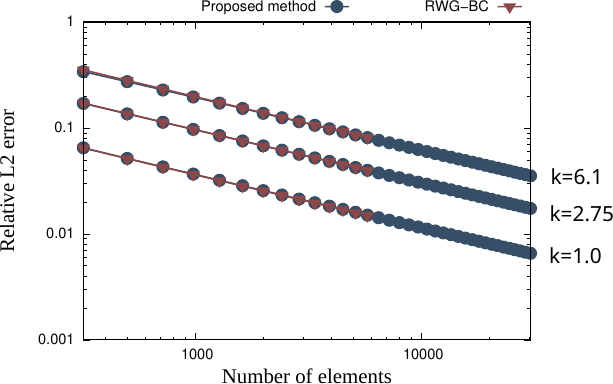}
    \caption{Relative errors of the numerical methods solving \eqref{eq:disc_reg_cfie_new} (blue) and \eqref{eq:disc_reg_cfie_conv} (red).
      The horizontal axis corresponds to the numbers of triangular meshes and the vertical axis to the $L_2$ relative error in \eqref{eq:relative_error}.
      Note that the numerical results for the conventional regularised CFIE in \eqref{eq:disc_reg_cfie_conv} are truncated up to $5780$ triangles
      since memory requirement for \eqref{eq:disc_reg_cfie_conv} is much more than that for \eqref{eq:disc_reg_cfie_new} due to the use of the BC function.}
    \label{figure:compare_rwg_bc}
  \end{center}
\end{figure}

We consider a scattering problem of a perfectly electric conductor~(PEC) with its radius $r=1$, which is illuminated with the plane wave.
The generalised minimul residual method~(GMRES) \cite{saad2003} with tolerance error $10^{-5}$ is used for solving linear equations.
For making implementaion simple, all the matrices in the discretised integral equations are directly computed, namely, all the elements of the matrices corresponding to discretisation of integral operators ($T$, $K$, $R$, and their tilde counterparts) are explicitly computed and fast methods such as the fast multipole method~(FMM) is not used.
We first compare the two regularised CFIE based on \eqref{eq:disc_reg_cfie_new} and \eqref{eq:disc_reg_cfie_conv}.
Fig.~\ref{figure:compare_rwg_bc} shows the relative error
\begin{align}\label{eq:relative_error}
        \frac{\|\bj_{\rm num} - \bj_{\rm ana}\|_{L^2_T}}{\|\bj_{\rm ana}\|_{L^2_T}}
\end{align}
of these numerical methods for several triangular meshes
where $\bj_{\rm num}$ is the solution obtained with the numerical methods and $\bj_{\rm ana}$ stands for the analytic solution with the Mie series.
Note that the numerical method based on equation \eqref{eq:disc_reg_cfie_conv} can be executed up to $5000$ triangles
since this method requires much more computational cost than that based on equation \eqref{eq:disc_reg_cfie_new}.
Table~\ref{table:rel_err} shows the relative errors and iteration nubmers of the GMRES with the number of triangles in a mesh $N=5120$ for the both numerical method and $N=32000$ for the method based on equation \eqref{eq:disc_reg_cfie_new}.
Fig.~\ref{figure:compare_rwg_bc} and Table~\ref{table:rel_err} imply that these two numerical methods show almost same accuracy for a mesh with same numbers of elements,
namely the computational cost of the proposed method for a fixed mesh is less than that of the conventional CFIE.
Also the iteration numbers of the GMRES for equation \eqref{eq:disc_reg_cfie_new} are almost same for the numerical result of different mesh size $N=5120$ and $N=32000$.
This indicates that the regularising operator $\mathcal{R}=\mathcal{T}^S_0$ works well in the proposed method.

\begin{figure}[htbp]
  \begin{center}
    \includegraphics[width=90mm]{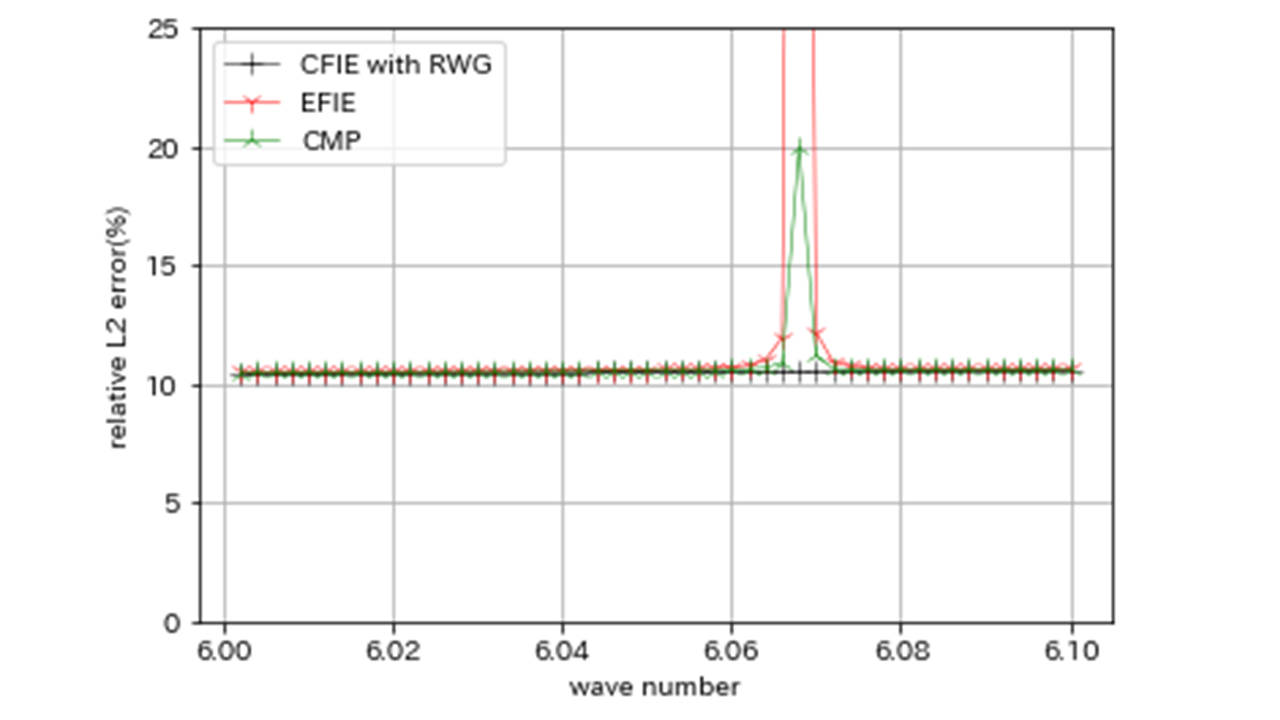}
    \caption{Relative error of the numerical methods.}\label{figure:comp_EC_error}
  \end{center}
\end{figure}

\begin{figure}[htbp]
  \begin{center}
    \includegraphics[width=90mm]{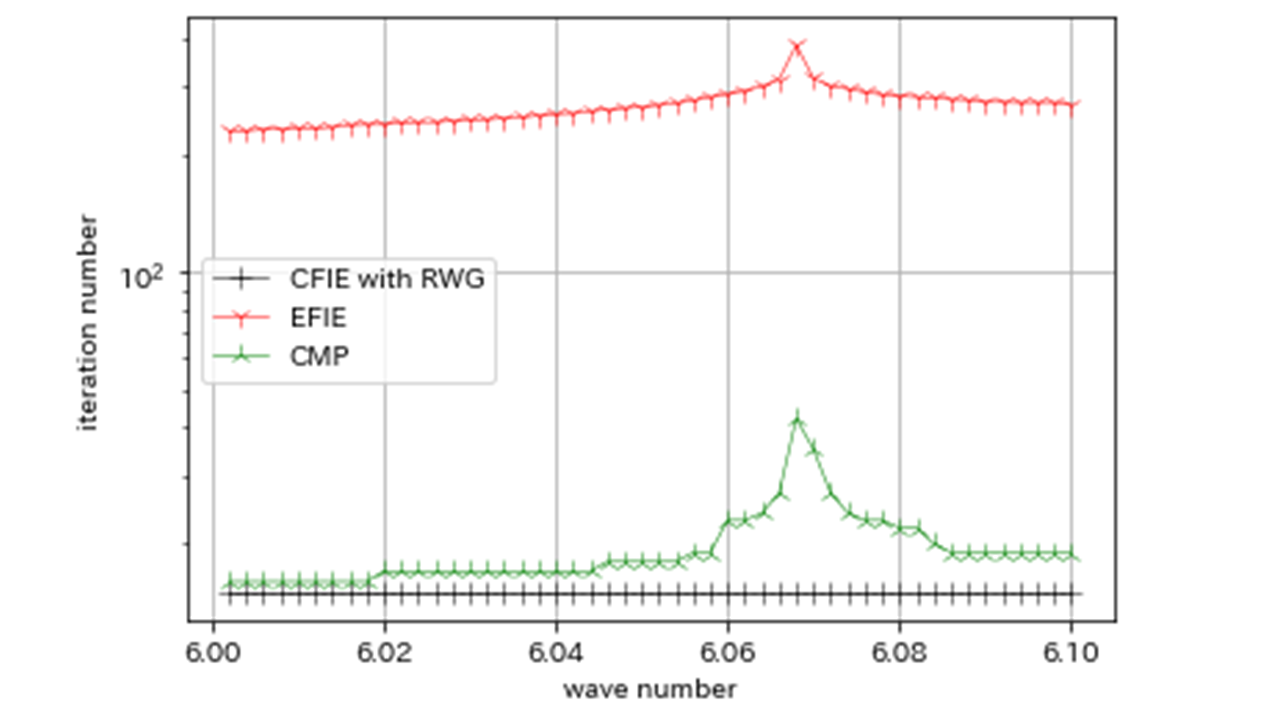}
    \caption{Iteration numbers of the GMRES}\label{figure:comp_EC_iteration}
  \end{center}
\end{figure}

Next we compare three numerical methods, the proposed CFIE solving equation \eqref{eq:disc_reg_cfie_new}, the standard EFIE and the EFIE preconditioned with the Calder\'{o}n multiplicative preconditioner~(CMP) \cite{andriulli2008}.
Fig.~\ref{figure:comp_EC_error} shows the relative errors \eqref{eq:relative_error} of the three numerical method for fixed mesh size $N=3380$.
In this problem the wave number $k\approx 6.07$ corresponds to one of the fictitious frequencies.
The relative error of the proposed method is stabl ein the range $6.0 \leq k \leq 6.1$ while those of the other two methods increase around the fictitious frequency.
Fig.~\ref{figure:comp_EC_iteration} shows the iteration numbers for the GMRES in the same setting as Fig.~\ref{figure:comp_EC_error}.
From these result we can see that the proposed regularised CFIE can successfully avoid the fictitious frequency
and is well-conditioned even around the fictitious frequency.


\section{Conclusion}
We proposed a discretisation method based on the Galerkin method for a regularised CFIE, which uses only the RWG basis function.
By using $\mathcal{R} = \mathcal{T}^S_0$ as the regularising operator, each operator in the regularised CFIE can be discretised without the dual basis function.
We also confirmed through a numerical example that the proposed method show good accuracy, and can avoid the ficticious-frequency problem.
In this article we focusd on proposing the formulation of the regularised CFIE and made simple implementation with the direct computation of matrices.
Hence the computational time of the numerical method was not discussed in section \ref{sec:numerical_example}.
Applying fast methods such as FMM to the proposed method and comparing the computatinal time of the proposed and conventional methods is a future plan.






%

\end{document}